\documentclass[twosided,12pt]{amsart}
\usepackage{amsmath,amsthm,amsfonts,amssymb}
\usepackage[pdftex,pdfborder={0 0 0}]{hyperref}
\usepackage{tikz}
\usepackage{float}
\usepackage{enumerate}
\usepackage{diagbox}
\usepackage{comment}
\usetikzlibrary{arrows,decorations.markings}
\usepackage{soul}
\usepackage{mathtools}
\usepackage{tikz}
\usetikzlibrary{graphs, quotes}
\usetikzlibrary{graphs.standard}

\tikzset{vertex/.style={circle,fill=black}}
\usetikzlibrary{spy}
\usetikzlibrary{fit,calc,positioning}
\usetikzlibrary{arrows.meta,patterns}
\usetikzlibrary{shapes.geometric}
\usetikzlibrary{decorations.pathreplacing,angles,quotes,calligraphy}
\usetikzlibrary{matrix,shapes,arrows,positioning}

\usepackage[skip=1pt plus1pt, indent=8pt]{parskip}

\newtheorem{theorem}{Theorem}[section]

\newtheorem{proposition}[theorem]{Proposition}
\newtheorem{corollary}[theorem]{Corollary}

\theoremstyle{definition}

\newtheorem{example}[theorem]{Example}

\newtheorem{remark}[theorem]{Remark}

\tikzstyle{arrow} = [thick,->,>=stealth]

\newlength{\Oldarrayrulewidth}

\newcommand{\Mod}[1]{\,(\mathrm{mod}\ #1)}

\def\m@th{\mathsurround=0pt}
\def\sm#1{\null\,\vcenter{\baselineskip9pt\lineskip.23ex\m@th
    \ialign{\hfil$\scriptstyle##$\hfil&&\ \hfil$\scriptstyle##$\hfil\crcr
    \mathstrut\crcr\noalign{\kern-\baselineskip}
    #1\crcr\mathstrut\crcr\noalign{\kern-\baselineskip}}}\,}
\def\smnp#1{\null\,\vcenter{\baselineskip9pt\lineskip.23ex\m@th
    \ialign{\hfil$\scriptstyle##$\hfil&&\ \ \hfil$\scriptstyle##$\hfil\crcr
    \mathstrut\crcr\noalign{\kern-\baselineskip}
    #1\crcr\mathstrut\crcr\noalign{\kern-\baselineskip}}}\,}

\tikzstyle{ipe stylesheet} = [
  ipe import,
  even odd rule,
  line join=round,
  line cap=butt,
  ipe pen normal/.style={line width=0.4},
  ipe pen heavier/.style={line width=0.8},
  ipe pen fat/.style={line width=1.2},
  ipe pen ultrafat/.style={line width=2},
  ipe pen normal,
  ipe mark normal/.style={ipe mark scale=3},
  ipe mark large/.style={ipe mark scale=5},
  ipe mark small/.style={ipe mark scale=2},
  ipe mark tiny/.style={ipe mark scale=1.1},
  ipe mark normal,
  /pgf/arrow keys/.cd,
  ipe arrow normal/.style={scale=7},
  ipe arrow large/.style={scale=10},
  ipe arrow small/.style={scale=5},
  ipe arrow tiny/.style={scale=3},
  ipe arrow normal,
  /tikz/.cd,
  ipe arrows, 
  <->/.tip = ipe normal,
  ipe dash normal/.style={dash pattern=},
  ipe dash dotted/.style={dash pattern=on 1bp off 3bp},
  ipe dash dashed/.style={dash pattern=on 4bp off 4bp},
  ipe dash dash dotted/.style={dash pattern=on 4bp off 2bp on 1bp off 2bp},
  ipe dash dash dot dotted/.style={dash pattern=on 4bp off 2bp on 1bp off 2bp on 1bp off 2bp},
  ipe dash normal,
  ipe node/.append style={font=\normalsize},
  ipe stretch normal/.style={ipe node stretch=1},
  ipe stretch normal,
  ipe opacity 10/.style={opacity=0.1},
  ipe opacity 30/.style={opacity=0.3},
  ipe opacity 50/.style={opacity=0.5},
  ipe opacity 75/.style={opacity=0.75},
  ipe opacity opaque/.style={opacity=1},
  ipe opacity opaque,
]
\definecolor{red}{rgb}{1,0,0}
\definecolor{blue}{rgb}{0,0,1}
\definecolor{green}{rgb}{0,1,0}
\definecolor{yellow}{rgb}{1,1,0}
\definecolor{orange}{rgb}{1,0.647,0}
\definecolor{gold}{rgb}{1,0.843,0}
\definecolor{purple}{rgb}{0.627,0.125,0.941}
\definecolor{gray}{rgb}{0.745,0.745,0.745}
\definecolor{brown}{rgb}{0.647,0.165,0.165}
\definecolor{navy}{rgb}{0,0,0.502}
\definecolor{pink}{rgb}{1,0.753,0.796}
\definecolor{seagreen}{rgb}{0.18,0.545,0.341}
\definecolor{turquoise}{rgb}{0.251,0.878,0.816}
\definecolor{violet}{rgb}{0.933,0.51,0.933}
\definecolor{darkblue}{rgb}{0,0,0.545}
\definecolor{darkcyan}{rgb}{0,0.545,0.545}
\definecolor{darkgray}{rgb}{0.663,0.663,0.663}
\definecolor{darkgreen}{rgb}{0,0.392,0}
\definecolor{darkmagenta}{rgb}{0.545,0,0.545}
\definecolor{darkorange}{rgb}{1,0.549,0}
\definecolor{darkred}{rgb}{0.545,0,0}
\definecolor{lightblue}{rgb}{0.678,0.847,0.902}
\definecolor{lightcyan}{rgb}{0.878,1,1}
\definecolor{lightgray}{rgb}{0.827,0.827,0.827}
\definecolor{lightgreen}{rgb}{0.565,0.933,0.565}
\definecolor{lightyellow}{rgb}{1,1,0.878}
\definecolor{black}{rgb}{0,0,0}
\definecolor{white}{rgb}{1,1,1}

\begin{document}

\title{Properties of Sub-Add Move Graphs}
 \author[\;\;\;\; Cesarz]{Patrick Cesarz}
\address{Patrick Cesarz, University of Wyoming, Laramie, WY, USA}
\email{pcesarz@uwyo.edu}
\author[Fiorini]{Eugene Fiorini}
\address{Eugene Fiorini, Rutgers University, DIMACS, Piscataway, NJ, USA}
\email{efiorini@dimacs.rutgers.edu}
\author[Gong]{Charles Gong}
\address{Charles Gong, Carnegie Mellon University, Pittsburgh, PA, USA}
\email{clgong@andrew.cmu.edu}
\author[Kelley]{Kyle Kelley}
\address{Kyle Kelley, Kenyon College, Gambier, OH, USA}
\email{kyleakelley@pm.me}
\author[Thomas]{Philip Thomas}
\address{Philip Thomas, Kutztown University, Kutztown, PA, USA}
\email{philipdthomas1@gmail.com}
\author[Woldar]{Andrew Woldar}
\address{Andrew Woldar, Villanova University, Villanova, PA, USA}
\email{andrew.woldar@villanova.edu}

 
\thanks{This research was supported in part by the National Science Foundation, grant MPS-2150299.}

\maketitle

\begin{abstract}
We introduce the notion of a move graph, that is,  a directed graph whose vertex set is a $\mathbb Z$-module $\mathbb Z_n^m$, and whose arc set is uniquely determined by the action $M\!:\!\mathbb Z_n^m\to \mathbb Z_n^m$ where $M$ is an $m\times m$ matrix with integer entries. We study the manner in which properties of move graphs differ when one varies the choice of cyclic group $\mathbb Z_n$. Our principal focus is on a special family of such graphs, which we refer to as ``sub-add move graphs.''
  \end{abstract}
 
\section{Introduction}

The origin of this research stems from a 2023 summer REU (Research Experience for Undergraduates) held at Moravian University in Bethlehem, PA. Initial investigations were loosely based on the notions of pebbling and chip-firing, two areas that have fairly robust histories replete with many diverse applications, e.g.\ see \cite{AS19,BCF03,BZ18,h1,Kl18}.
 
Early on, we felt it would be of interest to deviate from the standard rules of pebbling, whereby pebbles are transferred from a given vertex to some subset of its neighbors while exacting a ``fee'' for the transfer. It occurred to us that this problem could be of greater interest if the number of pebbles transferred would somehow depend on the particular vertices involved in the transaction. This eventually led to incorporating module structure into the problem so that pebble transfers could be dictated by some well-defined action of a group on the vertex set. 
 
We say a group $G$ acts on a directed graph $\Gamma$ if the elements of $G$ permute the vertices of $\Gamma$ while preserving arcs and non-arcs. 
 
  Over the past few decades, there has been extensive literature focused on various aspects of groups acting on graphs (e.g.\ geometric group actions on trees \cite{For02, LP97}, finite graphs \cite{Cor11, Mar98}, valency and girth \cite{Mar98, MN09}, and topological graphs \cite{DKQ11}, among many others). See \cite{AFM20, Con03, GRM97, KT18, SH18, SY17} for recent applications of groups acting on graphs to physiology, neural networks, chemistry, cryptography, and computational methods.

The notation and terminology we shall adopt are largely standard, e.g.\ see \cite{Bi93,GC01,Rotman} as excellent sources of background material.
  As is customary, we denote by 
 $V(\Gamma)$ and $A(\Gamma)$ the vertex set and arc set of a directed graph $\Gamma$, respectively. We define 
 the in-degree $\deg^-(v)$ and out-degree $\deg^+(v)$ of a vertex $v$ as follows:
\[\deg^-(v)=|\{w\in V(\Gamma)\mid (w,v)\in A(\Gamma)\}|,\]
\[\deg^+(v)=|\{w\in V(\Gamma)\mid (v,w)\in A(\Gamma)\}|.\]
Given $(v,w)\in A(\Gamma)$, we refer to $v$ as a \emph{parent} of $w$ and to $w$ as a \emph{child} of $v$. 
Thus, a vertex $v$ has $\deg^-(v)$ parents and $\deg^+(v)$ children.  

Given a directed graph $\Gamma$, we denote its underlying (undirected) graph by $\Gamma^u$.  A \emph{directed cycle} $C$ in $\Gamma$ is a cycle in 
$\Gamma^u$ with the property that in the subgraph of $\Gamma$ induced on $C$, one has $\deg^-(v)=\deg^+(v)=1$ for every vertex $v$ on $C$.  For notational convenience, we shall interpret $({v,v})\in A(\Gamma)$ as a directed $1$-cycle and a pair $({v,w}),({w,v})\in A(\Gamma)$ as a directed $2$-cycle.
 
 We say a directed graph $\Gamma$ is \emph{weakly connected} if $\Gamma^u$ is connected.  Similarly, a weakly connected component of $\Gamma$ corresponds to a connected component of $\Gamma^u$.

For us,  $\mathbb N$ will always denote the set of positive integers.  We adopt the notation ${\rm Mat}_m(\mathbb Z)$ for the algebra of $m\times m$ matrices with integer entries, as well as 
 $GL_m(\mathbb Q)$ for the group of  $m\times m$ invertible matrices with rational entries.  The order of a group  element $g$ will be designated  by ${\rm ord}(g)$. Lastly, ${\bf x}^T$ will denote the transpose of a vector $\bf x$. 

The balance of the paper is organized as follows. 
In Section 
\ref{sec:prelims}, we set the stage   by defining the notion of a move matrix and   its corresponding move graph.  In Section \ref{sec:3higherdim}, we investigate move graphs where the choice of move matrix is arbitrary. 
The sub-add move graph is introduced in Section \ref{sec:AddSub}.   
In Sections 
\ref{sec:n=2^r}--\ref{sec:n=p prime}, 
we treat special cases of the sub-add move graph. Specifically,  we  assume the underlying cyclic group $\mathbb Z_n$ has order $n=2^r$ (Sec.\ \ref{sec:n=2^r}), odd order $n$ (Sec.\ \ref{sec:n odd}), and odd prime order $p$ (Sec.\ \ref{sec:n=p prime}).

\section{Preliminaries}\label{sec:prelims}

Given $M\in {\rm Mat}_m(\mathbb Z)$, we consider the $m$-dimensional $\mathbb Z$-module $\mathbb Z_n^m$ where $\mathbb Z_n$ is the group of integers modulo $n$. If there exists a least positive integer $k$ for which $M^k$ acts as the identity on $\mathbb Z_n^m$, we refer to $k$ as the
\emph{$\mathbb Z_n$-order} of $M$.  
 This $\mathbb Z_n$-order plays a major role in the structure of $\Gamma_{M,\,n}$.

 
 We define the directed graph
  $\Gamma_{M,\,n}$ with vertex set $V(\Gamma_{M,\,n})$ and arc set $A(\Gamma_{M,\,n})$ as follows:
  \[V(\Gamma_{M,\,n}) = \{{\bf x}=(x_1,x_2,\dots,x_m) \mid x_i \in \mathbb{Z}_n\}\]
 \[A(\Gamma_{M,\,n}) = \{ ({\bf x}, {\bf y}) \mid {\bf x}, {\bf y}\in V(\Gamma_{M,\,n}), {\bf y}^T=M{\bf x}^T\}.\]
  In this context, we refer to $M$ as the \emph{move matrix} and to   
 $\Gamma_{M,\, n}$ as the \emph{move graph}.

We caution the reader that expressions of the form $ax$ where $a\in\mathbb Z$ and $x\in\mathbb Z_n$ should not be confused with the binary operation in  $\mathbb Z_n$ which is addition modulo $n$. 
Nevertheless, we can always interpret such expressions as elements of $\mathbb Z_n$. For example, $ax$ can be interpreted as the sum of $x$ with itself $a$ times if $a$ is positive, or $-x$ with itself $-a$ times  if $a$ is negative. 
 
  Our main objective is to study structural properties of $\Gamma_{M,\,n}$ and their dependence on  the move matrix $M$ and the cyclic group $\mathbb{Z}_n$.

\begin{example}\label{example:n=3}
Consider the move matrix $M$ given by 
\[M = \begin{pmatrix}
    0 & 0 & 1\\
    1 & 0 & 0\\
    0 & 1 & 0
\end{pmatrix}\]
As $M$ has order $3$ as an element of $GL_3(\mathbb Q)$, its $\mathbb Z_n$-order is $3$ for all $n\ge 2$.  We depict  the move graph $\Gamma_{M,\,3}$ in Figure \ref{fig:movegraphM3Z3}.
  
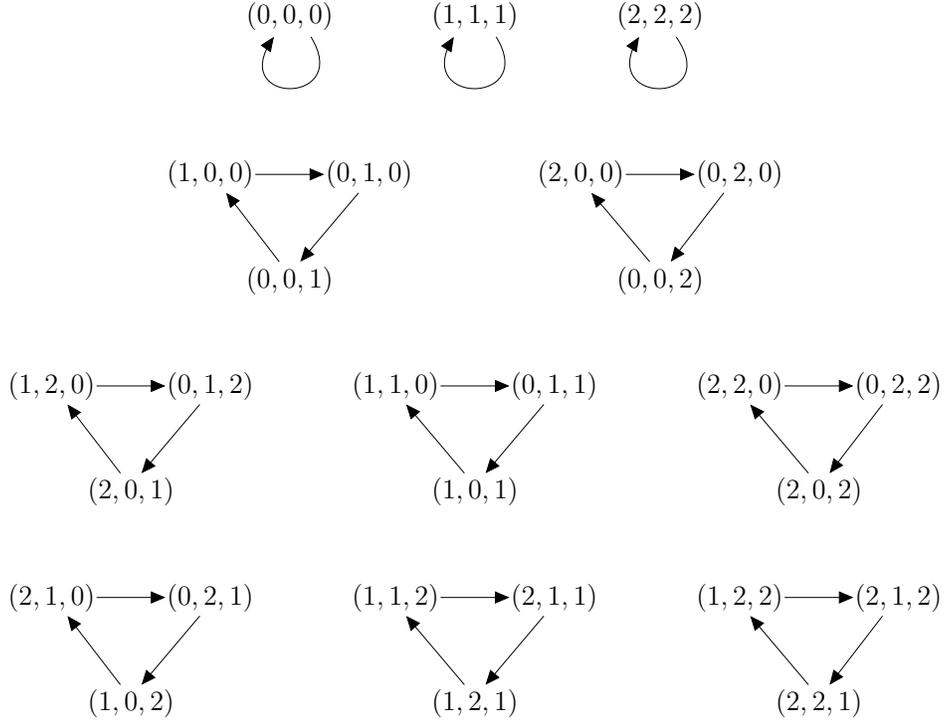
\begin{figure}[th]
\centering
\resizebox{1\textwidth}{!}{%
\begin{tikzpicture}[>=triangle 45,x=0.8cm,y=0.8cm]
    \node at (-3.5,7){$(0,0,0)$};
    \node at (0,7){$(1,1,1)$};
    \node at (3.5,7){$(2,2,2)$};
    \node at (-5,4){$(1,0,0)$};
    \node at (-2,4){$(0,1,0)$};
    \node at (2,4){$(2,0,0)$};
    \node at (5,4){$(0,2,0)$};
    \node at (-3.5,2){$(0,0,1)$};
    \node at (3.5,2){$(0,0,2)$};
    \node at (-8,0){$(1,2,0)$};
    \node at (-5,0){$(0,1,2)$};
    \node at (-1.5,0){$(1,1,0)$};
    \node at (1.5,0){$(0,1,1)$};
    \node at (5,0){$(2,2,0)$};
    \node at (8,0){$(0,2,2)$};
    \node at (-6.5,-2){$(2,0,1)$};
    \node at (0,-2){$(1,0,1)$};
    \node at (6.5,-2){$(2,0,2)$};
    \node at (-8,-4){$(2,1,0)$};
    \node at (-5,-4){$(0,2,1)$};
    \node at (-1.5,-4){$(1,1,2)$};
    \node at (1.5,-4){$(2,1,1)$};
    \node at (5,-4){$(1,2,2)$};
    \node at (8,-4){$(2,1,2)$};
    \node at (-6.5,-6){$(1,0,2)$};
    \node at (0,-6){$(1,2,1)$};
    \node at (6.5,-6){$(2,2,1)$};
    \draw[->](-4.15,4)--(-2.85,4);
    \draw[->](2.85,4)--(4.15,4);
    \draw[->](-2.2,3.65)--(-3.3,2.35);
    \draw[->](-3.7,2.35)--(-4.7,3.65);
    \draw[->](3.3,2.35)--(2.2,3.65);
    \draw[->](4.7,3.65)--(3.7,2.35);
    \draw[->](-7.15,0)--(-5.85,0);
    \draw[->](-.65,0)--(.65,0);
    \draw[->](5.85,0)--(7.15,0);
    \draw[->](-5.2,-.35)--(-6.3,-1.65);
    \draw[->](-6.7,-1.65)--(-7.7,-.35);
    \draw[->](1.3,-.35)--(.2,-1.65);
    \draw[->](-.2,-1.65)--(-1.3,-.35);
    \draw[->](6.3,-1.65)--(5.2,-.35);
    \draw[->](7.7,-.35)--(6.7,-1.65);
    \draw[->](-7.15,-4)--(-5.85,-4);
    \draw[->](-.65,-4)--(.65,-4);
    \draw[->](5.85,-4)--(7.15,-4);
    \draw[->](-5.2,-4.35)--(-6.3,-5.65);
    \draw[->](-6.7,-5.65)--(-7.7,-4.35);
    \draw[->](1.3,-4.35)--(.2,-5.65);
    \draw[->](-.2,-5.65)--(-1.3,-4.35);
    \draw[->](6.3,-5.65)--(5.2,-4.35);
    \draw[->](7.7,-4.35)--(6.7,-5.65);
    \draw[->] (-3.1,6.6) .. controls (-2.3,5.3) and (-4.7,5.3) .. (-3.8, 6.6);
    \draw[->] (0.4,6.6) .. controls (1.2,5.3) and (-1.2,5.3) .. (-0.4, 6.6);
    \draw[->] (3.8, 6.6)  .. controls (4.7,5.3) and (2.3,5.3) .. (3.1,6.6);
 
\end{tikzpicture}}
\vspace*{-3mm}
\caption{The move graph $\Gamma_{M,\,3}$ of Example \ref{example:n=3}.}\label{fig:movegraphM3Z3}
\end{figure}

\end{example}

\section{The general case}\label{sec:3higherdim}

\begin{theorem}\label{thm:Disjoint}
Let $M \in {\rm Mat}_{m}(\mathbb{Z})$  have $\mathbb Z_n$-order $k$. Suppose $M^{t}{\bf x}^T={\bf x}^T$ for some ${\bf x}\in V(\Gamma_{M,\,n})$, $t\in\mathbb N$. Then ${\bf x}$ is on a unique directed cycle in $\Gamma_{M,\,n}$ of length dividing $t$. In particular, all directed cycles in $\Gamma_{M,\,n}$ have length dividing $k$. 
\end{theorem}

\begin{proof}
We first prove every vertex of $\Gamma_{M,\,n}$ lies on a unique directed cycle. For this purpose, it suffices to show that ${\rm deg}^-({\bf x})={\rm deg}^+({\bf x})=1$ for all ${\bf x}\in V(\Gamma_{M,\,n})$. Of course, ${\rm deg}^+({\bf x})=1$ follows at once from the definition of $A(\Gamma_{M,\,n})$. 

The proof that ${\rm deg}^-({\bf x})=1$ for all ${\bf x}\in V(\Gamma_{M,\,n})$ hinges on the assumption that $M$ has finite $\mathbb Z_n$-order $k$. Here we observe that
$M^{k-1}$ is in ${\rm Mat}_{m}(\mathbb{Z})$  and this allows us to consider the move graph 
$\Gamma_{M^{k-1},n}$. By the above argument,
each vertex ${\bf x}\in V(\Gamma_{M^{k-1},n})$ has out-degree $1$,  which means for each ${\bf x}\in V(\Gamma_{M^{k-1},n})$ there exists a unique ${\bf y}\in V(\Gamma_{M^{k-1},n})$ for which
$({\bf x,y})\in A(\Gamma_{M^{k-1},n})$. But this implies   
${\bf y}^T=M^{k-1}{\bf x}^T$ and therefore 
$M{\bf y}^T=M^k{\bf x}^T={\bf x}^T$. It follows that 
$({\bf y,x})\in A(\Gamma_{M,\,n})$, i.e.\ ${\rm deg}^-({\bf x})= 1$ for all ${\bf x}\in V(\Gamma_{M,\,n})$. We  conclude that $\Gamma_{M,\,n}$ is a disjoint union of cycles.

Now suppose ${\bf x}$ lies on a directed $\ell$-cycle in $\Gamma_{M,\,n}$ for some $\ell\in\mathbb N$. 
Since cycles in $\Gamma_{M,\,n}$ are disjoint, it follows that $\ell$ is minimal such that
$M^\ell {\bf x}^T={\bf x}^T$.  Thus $\ell\le t$ whence by the Division Algorithm, 
$t=\ell q + r$ for some $q\in \mathbb{N}$ and $r\in\mathbb Z$ with $0 \le r < \ell$. 
But then
\[{\bf x}^T = M^t{\bf x}^T = M^{\ell q+r}{\bf x}^T = M^r((M^\ell)^q{\bf x}^T) = M^r{\bf x}^T.\]
Since $r\ne 0$ contradicts the minimality of $\ell$, we conclude that $r=0$, i.e.\ $\ell$ divides $t$ as claimed.
\end{proof}

\begin{proposition}\label{prop:LinearComb}
For fixed $\ell\in\mathbb N$, let $S$ be the set of all vertices on directed cycles in $\Gamma_{M,\,n}$ of length dividing $\ell$. Then any $\mathbb Z$-linear combination of vertices from $S$ yields a vertex on a cycle of length dividing $\ell$. Moreover, if $M$ has finite $\mathbb Z_n$-order, then for any ${\bf x} \in S$ and $s\in \mathbb Z$ with $\gcd(s,n)=1$, the vertices ${\bf x}$ and $s{\bf x}$ lie on cycles of equal length.
\end{proposition}

\begin{proof}
The first assertion is an immediate consequence of the linearity of $M$. To prove the second assertion, let ${\bf x}$ be on a cycle of length $\ell_1$ and let $s{\bf x}$ be on a cycle of length $\ell_2$. Then $M^{\ell_1}{\bf x}^T = {\bf x}^T$ and $sM^{\ell_2}{\bf x}^T=M^{\ell_2}(s{\bf x}^T) = s{\bf x}^T$.  Since $\gcd(s,n)=1$, the latter equality simplifies to $M^{\ell_2}{\bf x}^T = {\bf x}^T$.  
It therefore follows from Theorem \ref{thm:Disjoint} that
$\ell_1$ and $\ell_2$ divide one another, whence 
$\ell_1=\ell_2$ as claimed.
\end{proof}

\begin{theorem}\label{thm:subgraph}
    Let $n_1,n_2 \in \mathbb N$. Then $\Gamma_{M,\,n_1}$ embeds in $\Gamma_{M,\,n_1n_2}$ as an induced subgraph. 
\end{theorem}

\begin{proof}
   Define the mapping $f\!:\!V(\Gamma_{M,\,n_1}) \to V(\Gamma_{M,\,n_1n_2})$ by $f({\bf v})=n_2{\bf v}$.  We show $f$ is injective and preserves arcs and non-arcs.

    First suppose $f({\bf x}) = f({\bf y})$, i.e.\ $n_2{\bf x} = n_2{\bf y}$ in $\mathbb{Z}^m_{n_1n_2}$. Then ${\bf x} = {\bf y}$ in $\mathbb{Z}^m_{n_1}$ which proves $f$ is injective. Next observe that $n_2{\bf y}^T=(n_2{\bf y})^T=M(n_2{\bf x})^T=n_2M{\bf x}^T$  in $\mathbb{Z}^m_{n_1n_2}$ if and only if 
${\bf y}^T=M{\bf x}^T$  in $\mathbb{Z}^m_{n_1}$. This implies $({\bf x,y})\in  
A(\Gamma_{M,\,n_1})$ if and only if 
$(n_2{\bf x}, n_2{\bf y})\in  
A(\Gamma_{M,\,n_1n_2})$. The result follows.
\end{proof}
%
%

Let $\Delta_1$ and $\Delta_2$ be graphs with respective adjacency matrices $\mathcal A(\Delta_1)$ and $\mathcal A(\Delta_2)$. 
Recall that the {tensor product} $\Delta_1 \times \Delta_2$ is defined to be the graph 
with adjacency matrix $\mathcal A(\Delta_1 \times \Delta_2)$ equal to the Kronecker product $\mathcal A(\Delta_1) \otimes \mathcal A(\Delta_2)$.  


\begin{theorem}\label{thm:kbykKronecker}
Let $n_1,n_2\in \mathbb N$ with $\gcd(n_1,n_2) = 1$. Then the graph $\Gamma_{M,\,n_1n_2}$ is isomorphic to the tensor product $\Gamma_{M,\,n_1} \times \Gamma_{M,\,n_2}$.
\end{theorem}

\begin{proof}
Consider the mapping $f\!:\! V(\Gamma_{M,\,n_1}) \times V(\Gamma_{M,\,n_2}) \rightarrow V(\Gamma_{M,\,n_1n_2})$ defined by $f(\mathbf{x}, \mathbf{y}) = n_1\mathbf{y} + n_2\mathbf{x}$. We claim $f$ is a graph isomorphism. 

We first show $f$ is bijective. Suppose $f(\mathbf{x}, \mathbf{y}) = f(\mathbf{z}, \mathbf{w})$, i.e.\ $n_1\mathbf{y} + n_2 \mathbf{x} = n_1\mathbf{w} + n_2\mathbf{z}$. As $n_1n_2=0$ in 
$\mathbb Z_{n_1n_2}$, multiplying by $n_2$ yields $n_2^2 \mathbf{x} = n_2^2 \mathbf{z}$, which implies $n_2\mathbf{x} = n_2\mathbf{z}$ in $\mathbb{Z}^m_{n_1}$. Since $\gcd(n_1,n_2) = 1$, it follows that $\mathbf{x} = \mathbf{z}$ in $\mathbb{Z}^m_{n_1}$. A similar argument yields $\mathbf{y} = \mathbf{w}$ in $\mathbb{Z}^m_{n_2}$. Thus $(\mathbf{x},\mathbf{y}) = (\mathbf{z}, \mathbf{w})$, i.e.\ $f$ is injective. Since $|V(\Gamma_{M,\,n_1}) \times V(\Gamma_{M,\,n_2})| = n_1^mn_2^m=(n_1n_2)^m = |V(\Gamma_{M,\,n_1n_2})|$, we have that $f$ is surjective as well.

We next prove $f$ preserves arcs and non-arcs by establishing that $(n_1\mathbf{y} + n_2\mathbf{x}, n_1\mathbf{w}+n_2\mathbf{z}) \in A(\Gamma_{M,\,n_1n_2})$ if and only if $(\mathbf{x}, \mathbf{z}) \in A(\Gamma_{M,\,n_1})$ and $(\mathbf{y}, \mathbf{w}) \in A(\Gamma_{M,\,n_2})$.

 By definition, $(\mathbf{x} , \mathbf{z})  \in A(\Gamma_{M,\,n_1})$ implies $M\mathbf{x}^T = \mathbf{z}^T$ while $(\mathbf{y}, \mathbf{w})\in A(\Gamma_{M,\,n_2})$ implies $M\mathbf{y}^T = \mathbf{w}^T$. 
As a consequence, in $\mathbb{Z}^m_{n_1n_2}$ we have 
$M(n_2\mathbf{x})^T= n_2M\mathbf{x}^T = n_2\mathbf{z}^T$ and $M(n_1\mathbf{y})^T=n_1M\mathbf{y}^T = n_1\mathbf{w}^T$.
Summing these equalities yields  
$M(n_2\mathbf{x})^T + M(n_1\mathbf{y})^T = n_2\mathbf{z}^T + n_1\mathbf{w}^T$ or equivalently,  
$M(n_2\mathbf{x}+ n_1\mathbf{y})^T= 
(n_2\mathbf{z} + n_1\mathbf{w})^T$.  
We conclude that $(n_1\mathbf{y} + n_2\mathbf{x}, n_1\mathbf{w}+n_2\mathbf{z}) \in A(\Gamma_{M,\,n_1n_2})$.

 Conversely, suppose $(n_1\mathbf{y} + n_2\mathbf{x}, n_1\mathbf{w}+n_2\mathbf{z})\in A(\Gamma_{M,\,n_1n_2})$ in which case $n_1\mathbf{w}^T + n_2\mathbf{z}^T=M(n_1\mathbf{y}^T + n_2\mathbf{x}^T)=n_1M\mathbf{y}^T + n_2M\mathbf{x}^T$. 
Multiplying by $n_2$ gives $n_2^2\mathbf{z}^T=n_2^2M\mathbf{x}^T$ since $n_1n_2=0$ in $\mathbb{Z}_{n_1n_2}$. 
Thus $n_1n_2$ divides $n_2^2(\mathbf{z}^T-M\mathbf{x}^T)$ whence $n_1$ must divide $\mathbf{z}^T-M\mathbf{x}^T$ since $\gcd(n_1,n_2)=1$.
We conclude that $M\mathbf{x}^T=\mathbf{z}^T$ in $\mathbb Z_{n_1}^m$ whence $(\mathbf{x} , \mathbf{z})  \in A(\Gamma_{M,\,n_1})$.
By a symmetric argument, $(n_1\mathbf{y} + n_2\mathbf{x}, n_1\mathbf{w}+n_2\mathbf{z})\in A(\Gamma_{M,\,n_1n_2})$ implies $(\mathbf{y}, \mathbf{w})\in A(\Gamma_{M,\,n_2})$. It follows that $f$ preserves arcs and non-arcs, i.e.\ $f$ is a graph isomorphism as claimed.
\end{proof} 

\begin{theorem}\label{thm:GST}
 Let $M_1,M_2\in {\rm Mat}_m(\mathbb{Z})$ and $S\in {\rm Mat}_m(\mathbb Z)\cap GL_m(\mathbb Q)$ such that   $M_2=S^{-1}M_1S$.  Then for any $n\in\mathbb N$ with $\gcd(n, {\rm det}(S))=1$, the move graphs $\Gamma_{M_1,\,n}$ and $\Gamma_{M_2,\,n}$ are isomorphic. 
\end{theorem}

\begin{proof}
We define the map $f\!:\!V(\Gamma_{M_2,\,n}) \to V(\Gamma_{M_1,\,n})$ by $f({\bf v}) = {\bf v}S^T$. 
Note that ${\bf v}S^T\in V(\Gamma_{M_1,\,n})$ for all ${\bf v}\in V(\Gamma_{M_2,\,n})$ because $S^T\in 
{\rm Mat}_m(\mathbb{Z})$. We claim $f$ is a graph isomorphism.

Suppose first that $f({\bf x})={\bf 0}$, i.e.\ ${\bf x}S^T={\bf 0}$.  
Then $\sum_{i=1}^m x_i \boldsymbol{r}_i={\bf 0}$ where
$x_i$ is the $i$-th coordinate of ${\bf x}$, and $\boldsymbol{r}_i$ is the $i$-th row of $S^T$ taken modulo $n$. 
Assuming $\gcd(n,{\rm det}(S))=1$, it follows that the subset
$\{\boldsymbol{r}_1, \boldsymbol{r}_2,\dots, \boldsymbol{r}_m\}\subset \mathbb Z_n^m$ is linearly independent. Therefore $x_i=0$ for all $1\le i\le m$, i.e.\ ${\bf x}={\bf 0}$.


Now suppose $f({\bf v})=f({\bf w})$, i.e.\ ${\bf v}S^T=
{\bf w}S^T$.  Then $({\bf v}-{\bf w})S^T={\bf 0}$ whence
 it follows from above that
${\bf v}-{\bf w}={\bf 0}$, i.e.\ ${\bf v}={\bf w}$ in $\mathbb Z_n^m$. We thereby conclude that $f$ is injective.  Surjectivity, and therefore bijectivity,  now follow since $V(\Gamma_{M_1,\,n})=  
V(\Gamma_{M_2,\,n})$.

 Finally, observe that $({\bf x},{\bf y}) \in A(\Gamma_{M_2,\,n})$ if and only if ${\bf y}^T = M_2{\bf x}^T=S^{-1}M_1S{\bf x}^T$, i.e.\ if and only if    $(f({\bf y}))^T=S{\bf y}^T=M_1S{\bf x}^T=M_1(f({\bf x}))^T$. As this is equivalent to $(f({\bf x}), f({\bf y})) \in A(\Gamma_{M_1,\,n})$, we conclude that $f$ preserves arcs and non-arcs, hence  $f$ is a graph isomorphism. 
 \end{proof}

\begin{theorem}
Let $M \in {\rm Mat}_{m}(\mathbb{Z})$ and 
  $S\in {\rm Mat}_m(\mathbb Z)\cap GL_m(\mathbb Q)$ such that  $S^{-1}MS$ is the rational canonical form of $M$. Then for any prime $p$ with 
$\gcd({\rm det}(S),p)=1$ and for which $M$ has finite $\mathbb Z_p$-order $k$, the move graph $\Gamma_{M,\,p}$ contains a directed $k$-cycle. 
 \end{theorem}

\begin{proof}
Assume $p$ is a prime such that $\gcd({\rm det}(S),p)=1$ and suppose   
 $M$ has $\mathbb Z_p$-order $k$. 
By Theorem \ref{thm:GST}, it suffices to show 
$\Gamma_{S^{-1}MS,\,p}$ contains a directed $k$-cycle.   
%

Here $S^{-1}MS=M_1 \oplus M_2 \oplus \cdots \oplus M_f$ where for each $1\le i\le f$ and suitably chosen ${\bf v_i}\in \mathbb Z^m$, the $m_i\times m_i$ block matrix $M_i$ is expressed in terms of the basis  $\{{\bf v_i}^T,M{\bf v_i}^T,M^2{\bf v_i}^T,\dots, M^{m_i-1}{\bf v_i}^T\}$ for the $i$-th  $M$-invariant subspace of dimension $m_i$  over $\mathbb Q$.

 Now form the sum ${\bf v} = {\bf v_1} + {\bf v_2} + \cdots + {\bf v_f}$ where ${\bf v_i}$ is chosen in the above  manner for each $1\le i\le f$.  We 
 consider the corresponding sum taken modulo $p$, i.e.\ 
${\bf x} = {\bf x_1} + {\bf x_2} + \cdots + {\bf x_f}$ where 
${\bf x}, {\bf x}_i\in\mathbb Z_p^m$ with ${\bf x}\equiv {\bf v}\Mod{p}$ and ${\bf x_i}\equiv {\bf v_i}\Mod{p}$. 
By Theorem \ref{thm:Disjoint},  ${\bf x}$ lies on a directed $\ell$-cycle in $\Gamma_{S^{-1}MS,\,p}$ for some $\ell$ dividing $k$. This implies ${\bf x}^T = (S^{-1}MS)^\ell {\bf x}^T =(S^{-1}M^{\ell}S) {\bf x}^T = M_1^{\ell} {\bf x_1}^T \oplus  M_2^{\ell} {\bf x_2}^T \oplus \cdots \oplus M_f^{\ell} {\bf x_f}^T$ from which we deduce that $M_i^{\ell} {\bf x_i}^T = {\bf x_i}^T$ for each $1 \leq i \leq f$.   

Now let $k_i$ be the $\mathbb Z_p$-order of $M_i$, $1\le i\le f$. Since $S^{-1}MS$ and $M$ have the same $\mathbb Z_p$-order, it follows that $k={\rm lcm}(k_1,k_2, \dots,k_f)$. But $M_i^{\ell} {\bf x_i}^T = {\bf x_i}^T$ implies 
$k_i$ divides $\ell$ for each $1\le i\le f$. Thus $k$ divides $\ell$, whence $k=\ell$ from above. We conclude that  
$\Gamma_{S^{-1}MS,\,p}$ contains a directed $k$-cycle as desired.
\end{proof}

%
%

\section{The sub-add move
graph}\label{sec:AddSub}

Throughout this brief section, $M$ will denote  the so-called \emph{sub-add move matrix} given by 
\[M = \begin{pmatrix}  1 &\!\!  -1  \\ 1 &\;\,  1  \end{pmatrix}.\] 
The name derives from the action of $M$ on $\mathbb{Z}_n^2$, which results in subtraction in the first coordinate and addition in the second coordinate as indicated below.
\begin{equation*}
   \begin{pmatrix}  1 &\!\!  -1  \\ 1 &\;\,  1  \end{pmatrix}  \begin{pmatrix} a \\ b\end{pmatrix} = \begin{pmatrix} a-b \\ a+b \end{pmatrix}
\end{equation*}

Our main goal   is  to investigate the structure of the  move graph 
$\Gamma_{M,\,n}$ for certain values of $n$.  We refer to $\Gamma_{M,\,n}$
as a \emph{sub-add move graph}.

\begin{figure}[th]
\centering
\begin{tikzpicture}[>=triangle 45,x=1.0cm,y=.95cm]
    \node at (0,0){$(0,0)$};
    \node at (0,3){$(0,1)$};
    \node at (2.11,2.11){$(2,1)$};
    \node at (3,0){$(1,0)$};
    \node at (2.11,-2.11){$(1,1)$};
    \node at (0,-3){$(0,2)$};
    \node at (-2.11,-2.11){$(1,2)$};
    \node at (-3,0){$(2,0)$};
    \node at (-2.11,2.11){$(2,2)$};
    \draw[->](0.5,3)--(1.7,2.4);
    \draw[->](2.2,1.8)--(2.9,0.3);
    \draw[->](2.9,-0.3)--(2.2,-1.8);
    \draw[->](1.7,-2.4)--(0.5,-3);
    \draw[->](-0.5,-3)--(-1.7,-2.4);
    \draw[->](-2.2,-1.8)--(-2.9,-0.3);
    \draw[->](-2.9,0.3)--(-2.2,1.8);
    \draw[->](-1.7,2.4)--(-0.5,3);
    \draw[->](0.3,-0.3) .. controls (0.7,-1.2) and (-0.7,-1.2) .. (-0.3,-0.3);
\end{tikzpicture}
\caption{The sub-add move graph $\Gamma_{M,\,3}$.}\label{fig:movegraphZ3}
\end{figure}
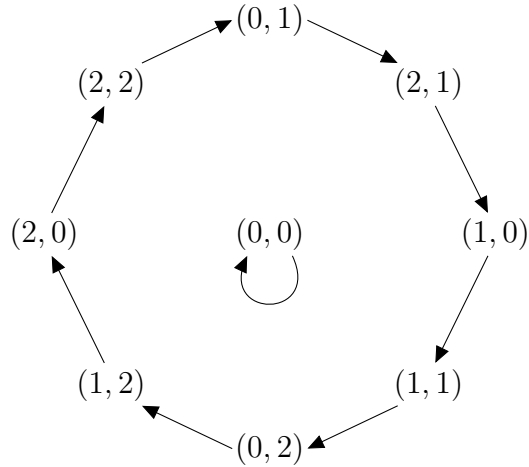

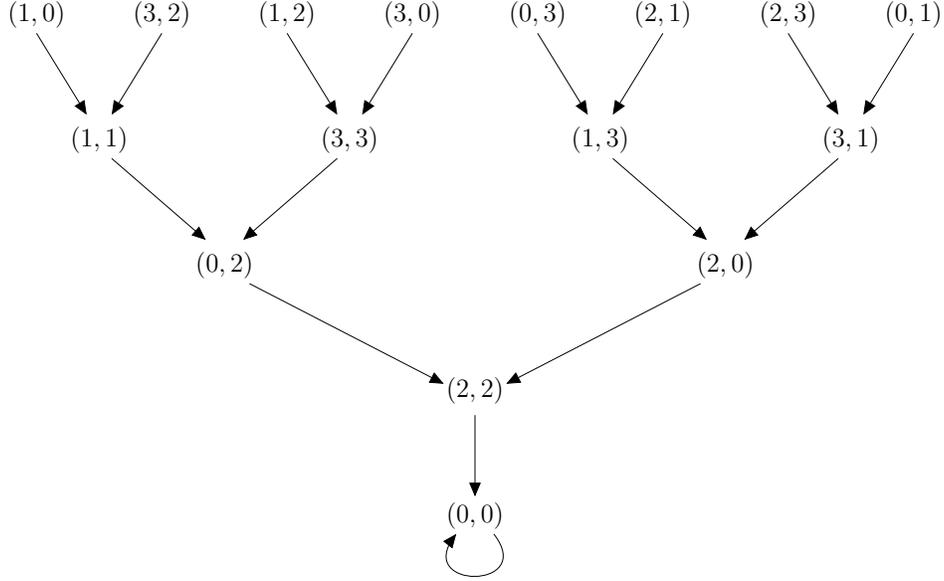
\begin{figure}[th]
\centering    
\resizebox{1\textwidth}{!}{%
\begin{tikzpicture}[>=triangle 45,x=1.0cm,y=1.0cm]
     \node at (-7,4){$(1,0)$};%
     \node at (-5,4){$(3,2)$};%
     \node at (-3,4){$(1,2)$};
     \node at (-1,4){$(3,0)$};
     \node at (1,4){$(0,3)$};
     \node at (3,4){$(2,1)$};
     \node at (5,4){$(2,3)$};
     \node at (7,4){$(0,1)$};
     \node at (-6,2){$(1,1)$};%
     \node at (-2,2){$(3,3)$};%
     \node at (2,2){$(1,3)$};%
     \node at (6,2){$(3,1)$};%
     \node at (-4,0){$(0,2)$};%
     \node at (4,0){$(2,0)$};%
     \node at (0,-2){$(2,2)$};%
     \node at (0,-4){$(0,0)$};%
     \draw[->](-7,3.7)--(-6.2,2.4);
     \draw[->](-5,3.7)--(-5.8,2.4);
     \draw[->](-3,3.7)--(-2.2,2.4);
     \draw[->](-1,3.7)--(-1.8,2.4);
     \draw[->](3,3.7)--(2.2,2.4);
     \draw[->](1,3.7)--(1.8,2.4);
     \draw[->](7,3.7)--(6.2,2.4);
     \draw[->](5,3.7)--(5.8,2.4);
     \draw[->](-5.8,1.7)--(-4.3,0.4);
     \draw[->](-2.2,1.7)--(-3.7,0.4);
     \draw[->](5.8,1.7)--(4.3,0.4);
     \draw[->](2.2,1.7)--(3.7,0.4);
     \draw[->](-3.6,-0.3)--(-0.5,-1.9);
     \draw[->](3.6,-0.3)--(0.5,-1.9);
     \draw[->](0,-2.4)--(0,-3.7);
     \draw[->](0.3,-4.3) .. controls (1,-5.2) and (-1,-5.2) .. (-0.3,-4.3);
\end{tikzpicture}}
\vspace*{-4mm}
\caption{The sub-add move graph $\Gamma_{M,\,4}$.}
 \label{fig:movegraphZ4}
\end{figure}

\begin{example}
In Figures \ref{fig:movegraphZ3} and \ref{fig:movegraphZ4} we depict the sub-add move graphs  $\Gamma_{M,\,3}$ and $\Gamma_{M,\,4}$, respectively.
\end{example}

\section{The sub-add move graph $\Gamma_{M,\,n}$, $n=2^r$}\label{sec:n=2^r}

 For fixed $r \in \mathbb N$, we may express the $\mathbb Z$-module $\mathbb Z_{2^r}^2$ in the following manner.
  \[\mathbb Z_{2^r}^2 = \{(2^tx,2^ty) \in \mathbb{Z}_{2^r}^2  \mid  \text{at least one of $x,y$ is odd}, 0\le t\le r\}.\]
This alternate description  makes it transparent that the subsets defined below form a partition of 
the vertex set $V(\Gamma_{M,\,2^r})=
\mathbb Z_{2^r}^2$.
 \[\begin{dcases*}
\;\;\mathcal P_{2t}= 
\{(2^t x, 2^t y)\mid \text{exactly one of $x$ and $y$ is odd}\}, & {$0\le t\le r-1$}\\[1ex]
\;\;\mathcal P_{2t+1}=\{2^{t}x, 2^{t} y)\mid \text{both $x$ and $y$ are odd}\}, & {$0\le t\le r-1$}\\[1ex]
\;\;\mathcal P_{2r} = \{(0,0)\}
\end{dcases*}\]
 
\begin{remark}
Technically, $x$ and $y$ are  elements of $\mathbb Z_{2^r}$ and not integers.  Nevertheless, we may regard them as integers when designating their parities. Indeed, since the group $\mathbb Z_{2^r}$ has even order, this designation is well-defined.  
\end{remark}

Occasionally, it will behoove us to  express subsets of the partition   simply as $\mathcal P_i$, $0\le i\le 2r$. This is purely for  notational convenience.  

%
 
%

%
%
%
%

 \begin{proposition}\label{prop:parentlevels}
     For every $0 \le i < 2r$, each vertex $(a_i,b_i) \in \mathcal{P}_i$ has a unique child $(a_{i+1},b_{i+1})$. Moreover, this child lies in   $\mathcal{P}_{i+1}$. For $i=2r$, the lone vertex $(0,0) \in \mathcal{P}_{2r}$ is its own unique child. 
 \end{proposition}

 \begin{proof}
Observe that $(2^{t}(x-y),2^{t}(x+y))$ is the unique child of $(2^{t}x,2^{t}y)\in\mathcal P_i$. Hence it remains to prove 
$(2^{t}(x-y),2^{t}(x+y))\in \mathcal P_{i+1}$. \\[.5mm]
Case 1:  $i=2t$, $0\le t<r$.  Since $(2^{t}x,2^{t}y)\in\mathcal P_{2t}$, $x$ and $y$ must have opposite parity. This means $x-y$ and $x+y$  are both odd, so it follows that   $(2^t(x-y),2^t(x+y)) \in \mathcal{P}_{2t+1}$.
 \\[.5mm]
 Case 2:     
   $i=2t+1$,  $0\le t<r$.  In this case, $(2^{t}x,2^{t}y)\in\mathcal P_{2t+1}$ so both $x$ and $y$ are odd.
Thus $x-y$ and $x+y$ are both even, so  the unique child 
$(2^{t}(x-y),2^{t}(x+y))$  is $(2^{t+1}(\frac{x-y}{2}), 2^{t+1}(\frac{x+y}{2}))$.  
      Now observe that $\frac{x-y}{2}+\frac{x+y}{2}=x$ which is odd. 
      This proves $\frac{x-y}{2}$ and $\frac{x+y}{2}$ must have opposite parity whence $(2^{t+1}(\frac{x-y}{2}), 2^{t+1}(\frac{x+y}{2}))\in \mathcal P_{2t+2}$ as desired.
   
 Finally it is obvious that the only child of $(0,0)$ is itself.   \end{proof}
 
\begin{proposition}\label{prop:parentless}
All vertices in $\mathcal P_0$ are parentless.
\end{proposition}

\begin{proof}
Let $(x,y)$ be an arbitrary vertex in $\mathcal P_0$ whence exactly one of $x, y$ is odd. This of course implies  $x+y$ must be odd.
Suppose now that $(a,b)$ is a parent of $(x,y)$.  
Then $(a-b,a+b)=(x,y)$ from which we obtain the contradiction $a=\frac{x+y}{2}$.  
\end{proof}

 Note that Propositions \ref{prop:parentlevels} and \ref{prop:parentless}
 establish the existence of levels in the sub-add move graph 
 $\Gamma_{M,\,2^r}$. That is to say, all vertices in $\mathcal P_0$ occur at the top level, their collective children occur at the next level,  and so on down to the bottom level,  which consists only of the single vertex $(0,0)$.

\begin{proposition}\label{prop:2parents}
Let $r>1$. Then $\big(2^{\lfloor \frac{i}{2} \rfloor}(x-y),2^{\lfloor \frac{i}{2} \rfloor}(x+y)\big) \in \mathcal{P}_{i+1}$ has exactly two parents for every $0 \leq i \le 2r-1$.


\end{proposition}

\begin{proof}
  Clearly, $(2^{\lfloor \frac{i}{2} \rfloor}x,2^{\lfloor \frac{i}{2} \rfloor}y),(2^{\lfloor \frac{i}{2} \rfloor}x+2^{r-1},2^{\lfloor \frac{i}{2} \rfloor}y+2^{r-1})$ are two distinct parents of $(2^{\lfloor \frac{i}{2} \rfloor}(x-y),2^{\lfloor \frac{i}{2} \rfloor}(x+y)) \in \mathcal{P}_{i+1}$.  
Let us assume this latter vertex has  a third parent. By Proposition \ref{prop:parentlevels}, this parent must lie in 
$\mathcal P_i$ and therefore be of the form  $(2^{\lfloor \frac{i}{2} \rfloor}x_1,2^{\lfloor \frac{i}{2} \rfloor}y_1)$.
  Since $\lfloor\frac{i}{2}\rfloor< r$, we thereby obtain $2^{\lfloor \frac{i}{2} \rfloor}(x_1-y_1) \equiv 2^{\lfloor \frac{i}{2} \rfloor}(x-y) \Mod{2^r}$ and $2^{\lfloor \frac{i}{2} \rfloor}(x_1+y_1) \equiv 2^{\lfloor \frac{i}{2} \rfloor}(x+y)\Mod{2^r}$. Since $r>1$, these together imply   
$2^{\lfloor \frac{i}{2} \rfloor}2x_1 \equiv 2^{\lfloor \frac{i}{2} \rfloor}2x\Mod{2^r}$ and $2^{\lfloor \frac{i}{2} \rfloor}2y_1 \equiv 2^{\lfloor \frac{i}{2} \rfloor}2y\Mod{2^r}$. 
  This gives rise to the following four possibilities:
    \begin{equation*}
        (2^{\lfloor \frac{i}{2} \rfloor}x_1,2^{\lfloor \frac{i}{2} \rfloor}y_1) = \begin{cases}
            (2^{\lfloor \frac{i}{2} \rfloor}x,2^{\lfloor \frac{i}{2} \rfloor}y) \\
            (2^{\lfloor \frac{i}{2} \rfloor}x + 2^{r-1},2^{\lfloor \frac{i}{2} \rfloor}y) \\
            (2^{\lfloor \frac{i}{2} \rfloor}x,2^{\lfloor \frac{i}{2} \rfloor}y + 2^{r-1}) \\
            (2^{\lfloor \frac{i}{2} \rfloor}x + 2^{r-1},2^{\lfloor \frac{i}{2} \rfloor}y + 2^{r-1})
        \end{cases}
    \end{equation*}
  But neither $(2^{\lfloor \frac{i}{2} \rfloor}x_1,2^{\lfloor \frac{i}{2} \rfloor}y_1)=(2^{\lfloor \frac{i}{2} \rfloor}x + 2^{r-1},2^{\lfloor \frac{i}{2} \rfloor}y)$ nor $(2^{\lfloor \frac{i}{2} \rfloor}x_1,2^{\lfloor \frac{i}{2} \rfloor}y_1)=(2^{\lfloor \frac{i}{2} \rfloor}x,2^{\lfloor \frac{i}{2} \rfloor}y + 2^{r-1})$ can occur because 
   $2^{\lfloor \frac{i}{2} \rfloor}x + 2^{\lfloor \frac{i}{2} \rfloor}y + 2^{r-1} \not\equiv 2^{\lfloor \frac{i}{2} \rfloor}x + 2^{\lfloor \frac{i}{2} \rfloor}y\Mod{2^r}$. 
As the remaining possibilities correspond to the  two aforementioned parents of 
$(2^{\lfloor \frac{i}{2} \rfloor}(x+y),2^{\lfloor \frac{i}{2} \rfloor}(y-x))$, the result follows when $i < 2r-1$.

   

The above argument also applies to the case $i=2r-1$, but here there is an interesting subtlety.  In this instance the child in question is   $(2^{r-1}(x-y),2^{r-1}(x+y))\in\mathcal P_{2r}=\{(0,0)\}$. From above, its only parents  are $(2^{r-1}x,2^{r-1}y)$ and 
$(2^{r-1}(x+1),2^{r-1}(y+1))$, however since $(2^{r-1}(x-y),2^{r-1}(x+y))=(0,0)$, we see that $x-y$ and $x+y$ must both be even implying that $x$ and $y$ have the same parity. If we assume $x$ and $y$ are both even, then $x+1$ and $y+1$ are both odd.  Alternatively, if we assume $x$ and $y$ are both odd, then $x+1$ and $y+1$ are both even.  In either case, one parent of $(0,0)$ is $(0,0)$ itself whereas the other parent lies in $\mathcal P_{2r-1}$. This completes the proof of the proposition.
\end{proof}

\begin{proposition}\label{prop:cardinality}
\color{black}{For $0 \leq i < 2r$, $|\mathcal{P}_i|= 2^{2r-i-1}$.}
\end{proposition}

\begin{proof}
As before, we consider individual cases based on the parity of $i$.\\[1mm]
Case 1: $i=2t$. 
Here there are $2^{r-t}$ choices for $x$, and  since $x$ and $y$ must have opposite parity there are $2^{r-t-1}$ independent choices for $y$. Thus, in total there are $2^{r-t} \cdot 2^{r-t-1} = 2^{2r-2t-1}$ choices for the pair $(x,y)$. Hence $|\mathcal{P}_{2t}| = 2^{2r-2t-1}=2^{2r-i-1}$. 
\\[1mm]
Case 2: $i=2t+1$.   
In this case $x$ and $y$ must both be odd. Thus there are $2^{r-t-1}$ choices for each of $x$ and $y$, hence 
$2^{r-t-1}  \cdot 2^{r-t-1} = 2^{2r-2t-2}$ pairs $(2^tx,2^ty)$ in total. Thus $|\mathcal{P}_{2t+1}| = 2^{2r-2t-2}= 2^{2r-i-1}$ in this case as well.
\end{proof}

Recall that a perfect binary tree $PBT_d$ is 
a type of binary tree where all non-leaf nodes have two children and all leaf nodes are at the same depth $d$. 
We refer to the unique parentless node $v$ of $PBT_d$ as its {root node}, and
we consider $PBT_d$ to be directed in the sense that there is a unique directed path from $v$ to each leaf node of $PBT_d$.   Finally, we call $PBT_d$ \emph{inverted} if the orientation of every arc of $PBT_d$ is reversed.


\begin{theorem}\label{thm:inverted PBT}
For each $r\in \mathbb N$, consider the subgraph $\Delta$ of 
$\Gamma_{M,\,2^r}$  induced on the vertex set $V(\Gamma_{M,\,2^r})\setminus \{(0,0)\}$. Then $\Delta$ is
an inverted $PBT_d$ of depth $d=2^r-1$. Moreover, 
one recovers $\Gamma_{M,\,2^r}$ from $\Delta$ by adding  the vertex $(0,0)$ along with the two arcs $((2^{r-1},2^{r-1}),(0,0))$ and $((0,0),(0,0))$.  
\end{theorem}

\begin{proof}
The structure of $\Delta$ follows directly from Propositions 
 \ref{prop:parentlevels}--\ref{prop:parentless}. 
To prove the latter assertion, we observe that 
the unique  childless vertex of $\Delta$ is $(2^{r-1},2^{r-1})$ and clearly $((2^{r-1},2^{r-1}),(0,0))\in A(\Gamma_{M,\,2^r})$.
\end{proof}

We once more refer the reader to Figure \ref{fig:movegraphZ4},  which illustrates  Theorem \ref{thm:inverted PBT} for the case $r=2$.

\section{The sub-add move graph $\Gamma_{M,\,n}$, $n$ odd}\label{sec:n odd}


\begin{theorem}\label{thm:n-odd}
    For $n\ge 3$ odd, $\Gamma_{M,\,n}$ is a disjoint union of directed cycles.
\end{theorem}

\begin{proof}
  By Theorem \ref{thm:Disjoint}, it suffices to show the move matrix $M$ has finite $\mathbb Z_n$-order.  However 
  $M^4=-4I$. 
  Thus $M$ has $\mathbb Z_n$-order dividing $4t$ where $t={\rm ord}(-4)$ in $\mathbb Z_n$.
%
%
%
 \end{proof}

\begin{figure}[th]
\centering
\resizebox{1\textwidth}{!}{%

\begin{tikzpicture}[>=triangle 45,x=1.0cm,y=1.0cm]

    \node at (-6,3.2){$(0,1)$};%
    \node at (-2,3.2){$(0,2)$};%
    \node at (6,3.2){$(0,4)$};%
    \node at (2,3.2){$(0,3)$};%
    
    \node at (-7,2){$(3,3)$};%
    \node at (-5,2){$(4,1)$};%
    \node at (-3,2){$(1,1)$};%
    \node at (-1,2){$(3,2)$};%
    \node at (3,2){$(2,3)$};%
    \node at (1,2){$(4,4)$};%
    \node at (5,2){$(2,2)$};%
    \node at (7,2){$(1,4)$};%
    
    \node at (-6,0.8){$(3,0)$};%
    \node at (-2,0.8){$(1,0)$};%
    \node at (6,0.8){$(2,0)$};%
    \node at (2,0.8){$(4,0)$};%
    
    \node at (-6,-0.8){$(1,3)$};%
    \node at (-2,-0.8){$(1,2)$};%
    \node at (2,-0.8){$(2,4)$};%
    
    \node at (6,-2){$(0,0)$};%
    \node at (-7,-2){$(2,1)$};%
    
    \node at (-5,-2){$(3,4)$};%
    
    \node at (-6,-3.2){$(4,2)$};%
    \node at (-2,-3.2){$(4,3)$};%
    \node at (2,-3.2){$(3,1)$};%
    
    \draw[->](-5.5,3)--(-5,2.3);
    \draw[->](-5,1.7)--(-5.5,1);
    \draw[->](-6.5,1)--(-7,1.7);
    \draw[->](-7,2.3)--(-6.5,3);
    
    \draw[->](-1.5,3)--(-1,2.3);
    \draw[->](-1,1.7)--(-1.5,1);
    \draw[->](-2.5,1)--(-3,1.7);
    \draw[->](-3,2.3)--(-2.5,3);
    
    \draw[->](2.5,3)--(3,2.3);
    \draw[->](3,1.7)--(2.5,1);
    \draw[->](1.5,1)--(1,1.7);
    \draw[->](1,2.3)--(1.5,3);

    \draw[->](6.5,3)--(7,2.3);
    \draw[->](7,1.7)--(6.5,1);
    \draw[->](5.5,1)--(5,1.7);
    \draw[->](5,2.3)--(5.5,3);

    \draw[->](-5.5,-1)--(-5,-1.7);
    \draw[->](-5,-2.3)--(-5.5,-3);
    \draw[->](-6.5,-3)--(-7,-2.3);
    \draw[->](-7,-1.7)--(-6.5,-1);

    \draw[<-](-2.5,-1) .. controls (-3,-2) ..(-2.5,-3);
    \draw[<-](-1.5,-3) .. controls (-1,-2) ..(-1.5,-1);

    \draw[<-](1.5,-1) .. controls (1,-2) ..(1.5,-3);
    \draw[<-](2.5,-3) .. controls (3,-2) ..(2.5,-1);

    \draw[->](6.3,-2.3) .. controls (7,-3.2) and (5,-3.2) .. (5.7,-2.3);
 
\end{tikzpicture}}
\caption{The sub-add move graph $\Gamma_{M,\,5}$ illustrates a special case of Theorem \ref{thm:n-odd}.}\label{fig:movegraphZ5}
\end{figure}

\begin{theorem}
    For $n\ge 3$ odd, the length of every  directed 
 cycle in $\Gamma_{M,\,n}$ divides $4\,\varphi(n)$ where $\varphi$ denotes Euler's totient function.
\end{theorem}

\begin{proof}
Since $M^4=-4 I$, we have $M^{4\varphi(n)}=(-4)^{\varphi(n)} I=4^{\varphi(n)} I$ since $\varphi(n)$ is even for all $n\ge 3$. 
Thus $M^{4\varphi(n)}({x},{y})^T=(4\varphi(n){x},4\varphi(n){y})^T$.  However, by Euler's totient theorem we have $4^{\varphi(n)}\equiv 1\Mod{n}$. Thus it follows that 
 $M^{4\varphi(n)}({x},{y})^T=({x},{y})^T$ whence
 the $\mathbb Z_n$-order of $M$ divides $4\,\varphi(n)$. The result now follows from Theorem \ref{thm:Disjoint}.
\end{proof}

\begin{theorem}\label{thm:(2L+1)2^k}
Let $n_1 \ge 3$ be odd and $n_2=2^k$. Then $\Gamma_{M,\,n_1n_2}$ contains $n_1^2$ vertex-disjoint copies of $\Gamma_{M,\,n_2}$. Furthermore, $\Gamma_{M,\,n_1n_2}$ and $\Gamma_{M,\,n_1}$ have the same number of weakly-connected components.
\end{theorem}

\begin{proof}
The first assertion is immediate from Theorem \ref{thm:kbykKronecker}. The second assertion follows since each directed cycle in $\Gamma_{M,\,n_1}$ extends to a unique weakly-connected component in $\Gamma_{M,\,n_1n_2}$.
\end{proof}

In Figure \ref{fig:movegraphZ6}, we provide an example of a move graph that illustrates Theorem \ref{thm:(2L+1)2^k}.

\begin{figure}[th]
\centering
\resizebox{1\textwidth}{!}{%
\begin{tikzpicture}[>=triangle 45,x=1.0cm,y=1.0cm]
    \node at (0,-1.5){$(0,0)$};
    \node at (0,0){$(1,1)$};
    \node at (-1.4,1.4){$(0,3)$};
    \node at (1.4,1.4){$(3,0)$};
    \draw[->](0.3,-1.8) .. controls (0.7,-2.7) and (-0.7,-2.7) .. (-0.3,-1.8);
    \draw[->](0,-0.3)--(0,-1.2);
    \draw[->](1.1,1.1)--(0.3,0.3);
    \draw[->](-1.1,1.1)--(-0.3,0.3);

    \node at (0,6){$(0,2)$};
    \node at (4.22,4.22){$(4,4)$};
    \node at (6,0){$(4,0)$};
    \node at (4.22,-4.22){$(2,4)$};
    \node at (0,-6){$(0,4)$};
    \node at (-4.22,-4.22){$(2,2)$};
    \node at (-6,0){$(2,0)$};
    \node at (-4.22,4.22){$(4,2)$};
    \draw[->](0.6,5.9)--(3.7,4.5);
    \draw[->](4.5,3.8)--(5.9,0.5);
    \draw[->](5.9,-0.5)--(4.4,-3.8);
    \draw[->](4,-4.6)--(0.7,-5.9);
    \draw[->](-0.6,-5.9)--(-3.7,-4.5);
    \draw[->](-4.5,-3.8)--(-5.9,-0.5);
    \draw[->](-5.9,0.5)--(-4.4,3.8);
    \draw[->](-4,4.6)--(-0.7,5.9);

    \node at (0,7.5){$(1,1)$};
    \node at (-1.4,8.9){$(1,0)$};
    \node at (1.4,8.9){$(4,3)$};
    \draw[->](0,7.2)--(0,6.3);
    \draw[->](1.1,8.6)--(0.3,7.8);
    \draw[->](-1.1,8.6)--(-0.3,7.8);

    \node at (5.62,5.62){$(1,3)$};
    \node at (5.62,7.3){$(2,1)$};
    \node at (7.7,5.62){$(5,4)$};
    \draw[->](5.3,5.3)--(4.6,4.6);
    \draw[->](5.62,7)--(5.62,6);
    \draw[->](7.2,5.62)--(6.15,5.62);

    \node at (8,0){$(5,1)$};
    \node at (9.4,1.4){$(0,1)$};
    \node at (9.4,-1.4){$(3,4)$};
    \draw[->](7.5,0)--(6.55,0);
    \draw[->](9.1,1.1)--(8.4,0.3);
    \draw[->](9.1,-1.1)--(8.4,-0.3);

    \node at (5.62,-5.62){$(3,1)$};
    \node at (5.62,-7.3){$(2,5)$};
    \node at (7.7,-5.62){$(5,2)$};
    \draw[->](5.3,-5.3)--(4.6,-4.6);
    \draw[->](5.62,-7)--(5.62,-6);
    \draw[->](7.2,-5.62)--(6.15,-5.62);

    \node at (0,-7.5){$(5,5)$};
    \node at (-1.4,-8.9){$(5,0)$};
    \node at (1.4,-8.9){$(2,3)$};
    \draw[->](0,-7.2)--(0,-6.3);
    \draw[->](1.1,-8.6)--(0.3,-7.8);
    \draw[->](-1.1,-8.6)--(-0.3,-7.8);

    \node at (-5.62,-5.62){$(5,3)$};
    \node at (-5.62,-7.3){$(4,5)$};
    \node at (-7.7,-5.62){$(1,2)$};
    \draw[->](-5.3,-5.3)--(-4.6,-4.6);
    \draw[->](-5.62,-7)--(-5.62,-6);
    \draw[->](-7.2,-5.62)--(-6.15,-5.62);

    \node at (-8,0){$(5,3)$};
    \node at (-9.4,1.4){$(4,5)$};
    \node at (-9.4,-1.4){$(1,2)$};
    \draw[->](-7.5,0)--(-6.55,0);
    \draw[->](-9.1,1.1)--(-8.4,0.3);
    \draw[->](-9.1,-1.1)--(-8.4,-0.3);

    \node at (-5.62,5.62){$(3,5)$};
    \node at (-5.62,7.3){$(4,1)$};
    \node at (-7.7,5.62){$(1,4)$};
    \draw[->](-5.3,5.3)--(-4.6,4.6);
    \draw[->](-5.62,7)--(-5.62,6);
    \draw[->](-7.2,5.62)--(-6.15,5.62);
    
\end{tikzpicture}}
\caption{The sub-add move graph $\Gamma_{M,\,6}$ illustrates a special case of Theorem \ref{thm:(2L+1)2^k}.} 
\label{fig:movegraphZ6}
\end{figure}

\section{The sub-add move graph $\Gamma_{M,\,p}$, $p$ an odd prime}\label{sec:n=p prime}

Throughout this section, $M$ will denote the sub-add move matrix and $p$ will be an odd prime. Let $t$ be the order of $-4$ in the multiplicative group of the finite field $GF(p)$. Since $M^4 = -4I$, it follows that $k=4t$ is the order of $M$. 

Our goal is  to determine the existence and length of  directed cycles in $\Gamma_{M,\,p}$ in terms of the eigenvalues of $M$, viz.\ $1+i$ and $1-i$ where $i^2=-1$. Without loss of generality, we may assume $s\le s'$ where $s={\rm ord}(1-i)$ and $s'={\rm ord}(1+i)$. Since $M^k=I$, it is clear that $s'=k$ and $s$ divides $k$.
Also, since $(1-i)^4 = -4$ and ${\rm ord}(-4) =t$, we have $s= t\gcd{(s, 4)}$. Clearly, this implies $s \in \{t, 2t, 4t\}$. 

\begin{proposition}\label{prop:eigenvalues}
The only possible cycle lengths in $\Gamma_{M,\,p}$
are $1$, $s$ and $k$, where $s$ and $k$ are not necessarily distinct.
\end{proposition}

\begin{proof}
Since $p$ is prime, we may regard $V(\Gamma_{M,\,p})=\mathbb Z_p^2$ as a $2$-dimensional vector space $V$ over  
$GF(p)$. 
Observe that $i\in GF(p)$ if and only if $-1$ is a square in $GF(p)$, i.e.\ if and only if $p\equiv 1\Mod{4}$.  
We treat this case first.

Let ${\bf v}^T$ and ${\bf w}^T$ be eigenvectors of $M$ corresponding to the eigenvalues $1+i$ and $1-i$, respectively. Since $1+ i, 1-i\in GF(p)$, 
we have that $\{{\bf v}^T, {\bf w}^T\}$ is a basis for 
$V$. 

Let now $C$ be an arbitrary directed cycle in $\Gamma_{M,\,p}$ and denote its length by $\ell$. Clearly, we may  assume $\ell\not\in \{1,k\}$. Let $\bf x$ be a vertex on $C$, in which case $M^\ell {\bf x}^T={\bf x}^T$. We may express ${\bf x}^T$ as $a{\bf v}^T + b{\bf w}^T$ for suitable scalars $a,b\in GF(p)$. Thus we obtain
{\small \[a{\bf v}^T+ b{\bf w}^T=M^\ell (a{\bf v}^T+ b{\bf w}^T)=
aM^\ell {\bf v}^T + b M^\ell {\bf w}^T=a(1+i)^\ell {\bf v}^T
+b(1-i)^\ell {\bf w}^T\]}

\noindent so that $a{\bf v}^T=a(1+i)^\ell {\bf v}^T$
and $b{\bf w}^T=b(1-i)^\ell {\bf w}^T$. 
The former equality implies either $a=0$ or $(1+i)^\ell = 1$. But $(1+i)^\ell = 1$ implies $\ell=k$, a case we have   ruled out by assumption.  Thus $a=0$, in which case
${\bf x}^T=b{\bf w}^T$, i.e.\ ${\bf x}^T$ is an eigenvector corresponding to the eigenvalue $1-i$.    
As such, ${\bf x}$ lies on a directed cycle of length $s$ since 
$M^s{\bf x}^T={\bf x}^T$. We thereby conclude that $\ell=s$ 
 by Theorem \ref{thm:n-odd}.  
 
The case $p\equiv 3\Mod{4}$ requires that  the above proof be modified slightly. Specifically,  
 $i$ is no longer an element of $GF(p)$ but lies in the unique quadratic extension $GF(p^2)$ of $GF(p)$. Here $i+1,i-1\in GF(p^2)$  are  algebraically conjugate over $GF(p)$, hence they have the same order $k$.  Note that in this case ${\bf v}^T, {\bf w}^T\notin V$, however $\{{\bf v}^T, {\bf w}^T\}$ still serves as a basis  for  a 2-dimensional vector space over $GF(p^2)$. As such, the linear combination ${\bf x}^T=a{\bf v}^T + b{\bf w}^T$ now allows  
 $a,b\in GF(p^2)$. Nevertheless, the  assumption that ${\bf x}\in V$ implies that $M^r {\bf x}^T\in V$ for all $r\in\mathbb N$. 
Thus we again reach $(1-i)^\ell =1$, from which we conclude that $\ell=k$.  
\end{proof}

We refer to a directed cycle of length $k$ in $\Gamma_{M,\,p}$ as a  \emph{primary cycle}. A \emph{secondary cycle} will be any directed cycle in $\Gamma_{M,\,p}$ that is neither primary nor a $1$-cycle. 

Note that primary cycles always exist. Indeed, if ${\bf v}^T$ is  
an eigenvector of $M$ corresponding to the eigenvalue $1+i$, then the vertex ${\bf v}$ will lie on a $k$-cycle since ${\rm ord}(1+i)=k$. In contrast, we see from   Proposition \ref{prop:eigenvalues} that secondary cycles exist  if and only if 
$s\in \{t,2t\}$ where $s={\rm ord}(1-i)$.  
 
%
%
%
%
%
%
%

Our next result elaborates on Proposition \ref{prop:eigenvalues} to a far greater extent.

\begin{theorem}\label{thm:secondary cycles}

\begin{enumerate}[\hspace*{9pt}\rm(a)]

    \item[]
    \item If $t$ is even, then $\Gamma_{M,\,p}$ has no secondary cycles. 
    
    \item If $t \equiv 1 \Mod{4}$ and $(1+i)^t=i$, then all secondary cycles have length $t$. 

    \item If $t \equiv 3 \Mod{4}$ and $(1+i)^t=i$, then all secondary cycles have length $2t$.

    \item If $t \equiv 1 \Mod{4}$ and $(1+i)^t = -i$, then all secondary cycles have length $2t$.

    \item If $t \equiv 3 \Mod{4}$ and $(1+i)^t=-i$, then all secondary cycles have length $t$.
    
\end{enumerate}
\end{theorem}

\begin{proof}
    Recall that ${\rm ord}(1+i)=4t$ which implies $(1+i)^t \in \{i,-i\}$. We first consider the case where $(1+i)^t=i$. As $1-i = -i(1+i)$ and $(1+i)^{3t}=-i$, we have $1-i=(1+i)^{3t+1}$. Observe that
\[{\rm ord}(1-i) = \frac{4t}{\gcd(3t+1,4t)} = \frac{4t}{\gcd(3t+1,t-1)} = \frac{4t}{\gcd(4,t-1)}.\]
 Thus, if $t$ is even we have $\gcd(4,t-1)=1$ whence ${\rm ord}(1-i)=4t$. We conclude that $\Gamma_{M,\,p}$ has no secondary cycles in this case.

    Next assume $t \equiv 1 \Mod{4}$. This implies $\gcd(4,t-1)=4$ whence ${\rm ord}(1-i)=t$. In this case,  all secondary cycles have length $t$. Similarly, if $t \equiv 3 \Mod{4}$, then $\gcd(4,t-1)=2$ which implies all secondary cycles have length $2t={\rm ord}(1-i)$.

      Finally, we consider the case $(1+i)^t=-i$. Here we have $1-i=(1+i)^{t+1}$ and therefore 
\[{\rm ord}(1-i) = \frac{4t}{\gcd(t+1,4t)} = \frac{4t}{\gcd(t+1,4)}.\]
 In this case, the roles of congruence are reversed in the sense that 
$t \equiv 1 \Mod{4}$ implies ${\rm ord}(1-i)=2t$ while $t \equiv 3 \Mod{4}$ implies ${\rm ord}(1-i)=t$. Accordingly, all secondary cycles in  $\Gamma_{M,\,p}$ have length $2t$ if $t \equiv 1 \Mod{4}$ or they have length $t$ if $t \equiv 3 \Mod{4}$. As above, if 
 $t$ is even then  $\gcd(t+1,4)=1$, in which case ${\rm ord}(1-i)=4t$, i.e.\ $\Gamma_{M,\,p}$ has no secondary cycles.  
\end{proof}


\begin{corollary}\label{cor:8dividesk}
  $\Gamma_{M,\,p}$ contains secondary cycles if and only if $8$ does not divide $k$. 
\end{corollary}

\begin{proof}
    This follows immediately from Theorem \ref{thm:secondary cycles} since $k$ is divisible by $8$ if and only if $t$ is even.
\end{proof}

\newpage

\begin{theorem}
\begin{enumerate}[\hspace*{9pt}\rm(a)]

    \item[]
        \item If $p \equiv 3 \Mod{8}$ or $p \equiv 7 \Mod{8}$, then secondary cycles do not exist in $\Gamma_{M,\,p}$.
        \item If $p \equiv 5 \Mod{8}$, then secondary cycles must exist in $\Gamma_{M,\,p}$.
    \end{enumerate}
\end{theorem}

\begin{proof}
 Both congruences in (a) imply $p\equiv 3\Mod{4}$, in which case $\Gamma_{M,\,p}$ contains  no secondary cycles by
 Proposition \ref{prop:eigenvalues}.
  Now suppose $p \equiv 5 \Mod{8}$. Clearly $8$ does not divide $p-1$, however  $k={\rm ord}(1+i)$ must divide $p-1$ since $1+i$ is an element of the multiplicative group of $GF(p)$.   
We conclude that $8$ does not divide $k$, whence it follows from Corollary \ref{cor:8dividesk} that  secondary cycles in $\Gamma_{M,\,p}$ must exist in this case. 
\end{proof}

\begin{remark}
At the present time, we have yet to discern any reasonable criteria that settle the question of existence  of secondary cycles  when 
 $p \equiv 1 \Mod{8}$. 
 Computer-generated data indicate that both outcomes are not only possible but occur with great frequency. Still, there doesn't appear to be any highly recognizable pattern in the data.
  \end{remark}


We close with the following result, which determines the number of cycles of varied length in the sub-add move graph $\Gamma_{M,\,p}$ when $p$ is an odd prime.  The proof is a straightforward counting argument so is left to the reader.

 \begin{proposition}
\label{prop:cyclecount}
\begin{enumerate}[\hspace*{9pt}\rm(a)]

    \item[]
        \item If $\,\Gamma_{M,\,p}$ has no secondary cycles, then there are $\frac{p^2 - 1}{k}$ primary cycles in $\Gamma_{M,\,p}$. 
        \item  If $\,\Gamma_{M,\,p}$ has secondary cycles, all such cycles have common length $s$ for some $s \in \{t, 2t\}$. Here there are $\frac{p-1}{s}$ secondary cycles and $\frac{p^2-p}{k}$ primary cycles in $\Gamma_{M,\,p}$.
    \end{enumerate}
\end{proposition}

\begin{remark}
    The entry A363894 in the On-Line Encyclopedia Integer Sequences \cite{oeis} reflects Proposition \ref{prop:cyclecount} in conjunction with Theorem \ref{thm:(2L+1)2^k}. Note that what is referred to as the ``halved addsub configuration graph'' in A363894 is the ``sub-add move graph'' in our present terminology.
\end{remark}  

While many of our results in this section do not extend to proper odd prime powers, a number of them do. We intend to pursue this line of investigation in a sequel.


\section*{Conflict of interest}
The authors declare that they have no conflict of interest.

\section*{Data Availability}
Data sharing is not applicable to this article as no datasets were generated or analyzed during the current study.

\section*{Funding information}

 This research was supported in part by the National Science Foundation, grant MPS-2150299. 
    
 \section*{Author contribution}
 
 Authors contributed to the study conception and design. Material preparation, data collection and analysis were performed by Patrick Cesarz, Charles Gong, Eugene Fiorini, Kyle Kelley, Philip Thomas, and Andrew Woldar. The first draft of the manuscript was written by Andrew Woldar and all authors commented on previous versions of the manuscript. All authors read and approved the final manuscript.

\bibliographystyle{plain}
\bibliography{bibliography}


\end{document}